\documentclass{svjour3}                     
\smartqed  
\usepackage{latexsym}
\usepackage{amssymb}
\usepackage{amsmath}
\usepackage{rotating}
\usepackage[misc]{ifsym}

\newtheorem{them}{Theorem}[section]
\newtheorem{lem}{Lemma}[section]

\newtheorem{coro}{Corollary}[section]
\newtheorem{pro}{Proposition}[section]

\journalname{Optim. Letters}
\begin{document}

\title{On Constraint Qualifications of a Nonconvex Inequality\thanks{This research was supported by the National Natural Science Foundations of P. R. China (grant 11401518) and the Fok Ying-Tung Education Foundation (grant 151101). The second author
was partially supported by the Grant MOST 105-2115-M-039-002-MY3.}}

\titlerunning{On Constraint Qualifications of a Nonconvex Inequality}



\author{Zhou Wei \and Jen-Chih Yao}


\institute{ Zhou Wei(\Letter) \at Department of Mathematics, Yunnan University, Kunming 650091, People's Republic of China\\ \email{wzhou@ynu.edu.cn} \and Jen-Chih Yao \at Center for General Education, China Medical University, Taichung 40402, Taiwan\\ \email{yaojc@mail.cmu.edu.tw}}

\date{Received: date / Accepted: date}

\maketitle

\begin{abstract}
In this paper, we study constraint qualifications for the nonconvex inequality defined by a proper lower semicontinuous function. These constraint qualifications involve the generalized construction of normal cones and subdifferentials. Several conditions for these constraint qualifications are also provided therein. When restricted to the convex inequality, these constraint qualifications reduce to basic constraint qualification (BCQ) and strong BCQ studied in [SIAM J. Optim., 14(2004), 757-772] and [Math. Oper. Res., 30 (2005), 956-965].


\end{abstract}

\section{Introduction}

Constraint qualifications (involving epigraph or subdifferential) have been widely studied and extensively used in various aspects of optimization and mathematical programming (see \cite{4a,6a,12a,13a,17a,22a,FLLY,FLN,FLY,FZ,LZH,Zh1,Zh2} and references therein). For example, constraint qualifications were used to study Fenchel duality and the formula of subdifferentials of convex functions (cf. \cite{4a,6a,12a,13a}). One Farkas-type constraint qualification has been proved to be useful in convex programming and DC (difference of two convex functions) programming (cf. \cite{17a,22a,LZH}). Constraint qualifications involving epigraphs were applied to the extended Farkas lemma and Lagrange duality in convex programming (cf. \cite{FLLY,FLN,FLY}). Constraint qualification were also applied to the study of optimality conditions in convex and DC optimization problems (cf. \cite{FLN,FZ}). Recently one type of closed cone constraint qualifications was used to study Lagrange duality in quasi-convex programming (cf. \cite{Zh1,Zh2}).

Since constraint qualifications are proved to play an important role in optimization and mathematical programming, there have been fruitful works on this topic and several types of constraint qualifications  have been extensively studied. One of the most important constraint qualifications in convex programming is basic constraint qualification (BCQ). BCQ is a fundamental concept in mathematical programming and has been widely discussed by many authors (cf. \cite{BB,De,HL,H,Li,LN,LN1,LN2,LNS,NY,Za,Z1,Z2} and references therein). For a continuous convex function $\phi$ defined on a Banach space $X$,  recall that convex inequality $ \phi(x)\leq 0$ is said to satisfy BCQ at $x\in S:=\{u\in X: \phi(u)\leq 0\}$ if
\begin{equation}\label{1.1}
  N(S,x)=[0,+\infty)\partial\phi(x)
\end{equation}
where $N(S,x)$ and $\partial\phi(x)$ refer to the normal cone and the subdifferential in the sense of convex analysis, respectively. It has been well recognized that BCQ closely relates to many important notions in convex analysis and optimization. Deutsch \cite{De} studied strong conical hull intersection property (strong CHIP) in best constraint approximation and showed that strong CHIP is a geometric version of BCQ. Li \cite{Li} studied BCQ and Abadie constraint qualification (ACQ) for differentiable convex inequalities and proved the equivalence between ACQ and BCQ. Subsequently Li, Nahak and Singer \cite{LNS} investigated BCQ and ACQ for semi-infinite systems of convex inequalities. For a convex semi-infinite optimization problem, Li, Ng and Pong \cite{LN2} showed that the BCQ is necessary and sufficient for a feasible point to be an optimal solution. Fang, Li and Ng \cite{FLN} proved that BCQ is equivalent to the KKT optimality conditions for the convex semi-infinite optimization problem.

In 1997 Lewis and Pang \cite{Le} studied metric regularity of a convex inequality and proved that BCQ is a necessary condition for metric regularity. For characterizing metric regularity, Zheng and Ng \cite{Z1} introduced and studied the concept of strong basic constraint qualification (strong BCQ) for a convex inequality. Recall from \cite{Z1} that convex inequality $\phi(x)\leq 0$ is said to satisfy strong BCQ at $x\in {\rm bd}(S)$ if there exists $\tau\in(0,+\infty)$ such that
\begin{equation}\label{1.2}
   N(S,x)\cap B_{X^*}\subset [0,\tau]\partial\phi(x)
\end{equation}
where ${\rm bd}(S)$ denotes the boundary of solution set $S$ and $B_{X^*}$ denotes the closed unit ball of dual space $X^*$. This constraint qualification, also named as $\tau$-strong BCQ, has been used to characterize metric regularity of convex inequality $\phi(x)\leq 0$ (see \cite[Theorem 2.2]{Z1} for more details). Further Hu \cite{H} studied the difference between strong BCQ in \eqref{1.2} and BCQ in \eqref{1.1} and introduced the concept of the end set to characterize this difference. It has been showed that  strong BCQ  is equivalent to BCQ and positive distance of the origin to the end set of subdifferential $\partial\phi(x)$ (cf. \cite[Theorem 3.1]{H}). Subsequently Zheng and Ng \cite{Z2} discussed BCQ and strong BCQ for a convex constraint system and applied strong BCQ to prove dual characterizations for metric subregularity of convex constraint system. Recently Huang and Wei \cite{HW} studied metric subregularity of convex constraint system by primal equivalent conditions and demonstrated that these primal conditions can characterize strong BCQ for this case (cf. \cite[Theorem 3.1 and Propoition 4.1]{HW}).

Note that nonconvex functions frequently appear in mathematical programming and several types of normal cones and subdifferentials for the nonconvex case have been extensively studied in optimization and its applications (cf. \cite{C,M,RW} for generalized constructions of normal cones and subdifferentials). Thus it is natural to further study constraint qualifications  by these normal cones and subdifferentials for a nonconvex inequality. Our main aim is to consider constraint qualifications involving normal cones nd subdifferentials for a nonconvex equality and extend results in \cite{H,Z1} to the nonconvex case. Motivated by this, we mainly study constraint qualifications involving Clarke/Fr\'echet normal cones and subdifferentials for a nonconvex inequality and aim to provide conditions for ensuring these constraint qualifications. It is proved that some results on BCQ and strong BCQ given in \cite{H,Z1} are also valid for the nonconvex case.

The  rest is the organization of this paper. In Section 2, we give some definitions and preliminaries used in our analysis. Our notation is basically standard and conventional in optimization and mathematical programming. Section 3 is devoted to constraint qualifications (that are named as Clarke BCQ and Clarke strong BCQ) for a nonconvex inequality. We first study Clarke BCQ and strong BCQ of the nonconvex inequality which is defined by a local Lipschitz function. Then we use Clarke singular subdifferential to consider the case where the inequality is defined by a lower semicontinuous function (not local Lipschitz necessarily). Several necessary and/or sufficient conditions for these constraint qualifications are also presented therein. In Section 4, we further discuss constraint qualifications for the nonconvex inequality by Fr\'echet normal cone and subdifferential.  The conclusion of this paper is presented in Section 5.

\setcounter{equation}{0}

\section{Preliminaries}
Let $X, Y$ be Banach spaces with the closed unit balls denoted by
$B_{X}$ and $B_Y$, and let $X^*, Y^*$ denote the dual spaces of $X$
and $Y$, respectively. For a subset $A$ of $X$, let $\overline
A$ and ${\rm bd}(A)$ denote the closure and the boundary of $A$, respectively. For any $x\in X$, we denote by $d(x, A):=\inf_{a\in A}\|x-a\|$ the distance of $x$ to $A$.

Let $A$ be a closed subset of $X$ and $a\in A$. We denote by $T(A,a)$ and $T_c(A,a)$ the Bouligand contingent cone and the Clarke tangent cone of $A$ at $a$, respectively which are defined by
$$
T(A,a):=\mathop{\rm Limsup}\limits_{t\rightarrow
0^+}\frac{A-a}{t}\ \ {\rm and} \ \ T_c(A,a):=\mathop{\rm Liminf}\limits_{x\stackrel{A}{\rightarrow}a,t\rightarrow
0^+}\frac{A-x}{t},
$$
where $x\stackrel{A}{\rightarrow}a$ means that $x\rightarrow a$ with
$x\in A$. Thus, $v\in T(A,a)$ if and only if there exist $t_n\rightarrow 0^+$ and $v_n\rightarrow v$ such that $a+t_nv_n\in A$ for all $n$, and $v\in T_c(A,a)$ if and only if for any
$a_n\stackrel{A}{\rightarrow}a$ and any $t_n\rightarrow 0^+$, there exists $v_n\rightarrow v$ such that $a_n+t_nv_n\in A$ for all $n$.

We denote by
$$N_c(A,a):=\{x^*\in X^*:\;\langle x^*,h\rangle\leq0\;\;\forall h\in
T_c(A,a)\}$$
the Clarke normal cone of $A$ at $a$ and by
$$
\hat N(A, a):=\left\{x^*\in X^*: \limsup_{x\stackrel{A}\rightarrow a}\frac{\langle x^*, x-a\rangle}{\|x-a\|}\leq 0\right\}
$$
the Fr\'echet normal cone of $A$ at $a$. It is easy to verify that $\hat N(A,a)$ is a norm-closed convex cone of $X^*$ and
$$
\hat N(A, a)\subset N_c(A,a).
$$
When $X$ is finite-dimensional, one has
\begin{equation}\label{2.1}
  \hat N(A,a)=(T(A,a))^{\circ}:=\{x^*\in X^*: \langle x^*,h\rangle\leq 0\ \forall h\in T(A,a)\}
\end{equation}
where $(T(A,a))^{\circ}$ is the nonnegative polar of $T(A,a)$. If $A$ is convex, both the Clarke normal cone and the Fr\'eceht normal cone reduce to that in the sense of convex analysis; that is,
$$
\hat N(A, a)=N_c(A,a)=N(A,a)=\{x^*\in X^*:\langle x^*,x-a\rangle\leq 0 \ \ \forall x\in A\}.
$$

Let  $\phi :X\rightarrow \mathbb{R}\cup\{+\infty\}$ be a proper
lower semicontinuous function and $x\in
\mathrm{dom}(\phi):=\{y\in X:\phi(y)<+\infty\}$. We denote by
$$
\hat\partial\phi(x):=\{x^*\in X^* : (x^*, -1)\in \hat N({\rm epi}(\phi), (x, \phi(x)))\}
$$
the Fr\'echet subdifferential of $\phi$ at $x$, where ${\rm epi}(\phi):=\{(u,r)\in X\times\mathbb{R}: \phi(u)\leq r\}
$ is the epigraph of $\phi$. It is known that $\hat\partial\phi(x)$ is a norm-closed convex (empty possibly) subset of $X^*$ and
$$
\hat\partial\phi(x)=\left\{x^*\in X^* : \liminf_{y\rightarrow x}\frac{\phi(y)-\phi(x)-\langle x^*, y-x\rangle}{\|y-x\|}\geq 0\right\}.
$$

We denote by $\partial_c\phi(x)$ and $\partial_c^{\infty}\phi(x)$ the Clarke subdifferential and the Clarke singular subdifferential of $\phi$ at $x$, respectively and they are defined as
\begin{eqnarray*}
\partial_c\phi(x):&=&\{x^*\in X^* : (x^*, -1)\in N_c({\rm epi}(\phi), (x, \phi(x)))\},\\
\partial_c^{\infty}\phi(x):&=&\{x^*\in X^* : (x^*, 0)\in N_c({\rm epi}(\phi), (x, \phi(x)))\}.
\end{eqnarray*}
It is known that
$$
\hat\partial\phi(x)\subset\partial_c\phi(x).
$$
When $\phi$ is convex, the Clarke subdifferential of $\phi$ at $x$ coincides with the Fr\'echet subdifferential and both reduce to that in the sense of convex analysis; that is,
\begin{eqnarray*}
\partial_c\phi(x)=\hat\partial\phi(x)=\partial\phi(x)=\{x^*\in X^* : \langle x^*, y-x \rangle\leq\phi(y)-\phi(x) \ \forall y\in X\},
\end{eqnarray*}
and
$$
\partial_c^{\infty}\phi(x)=\partial^{\infty}\phi(x)=\{x^*\in X^* : \langle x^*, y-x \rangle\leq 0\ \forall y\in{\rm dom}(\phi)\}.
$$

For any $h\in X$, we denote by $\phi^{\uparrow}(x,h)$ the generalized Rockafellar directional derivative of $\phi$
at $x$ along the direction $h$ which is
defined by (see \cite{C})
$$
\phi^{\uparrow}(x,h):=\lim\limits_{\varepsilon\downarrow
0}\limsup\limits_{z\stackrel{\phi}{\rightarrow}x,
\,t\downarrow0}\inf_{w\in h+\varepsilon
B_{X}}\frac{\phi(z+tw)-\phi(z)}{t},
$$
where  $ z\stackrel{\phi}\rightarrow x$ means that $z\rightarrow x$
and $\phi(z)\rightarrow \phi(x)$.

It is known from \cite{C} that
$$
\partial_c\phi(x)=\{x^*\in X^*:\langle x^*, h\rangle\leq \phi^{\uparrow}(x,h)\,\,\forall h\in
X\}.
$$
When $\phi$ is local Lipschitz around $x$, $\phi^{\uparrow}(x,h)$ reduces to the Clarke
directional derivative; that is
$$
\phi^{\uparrow}(x,h)=\phi^{\circ}(x,h):=\limsup\limits_{z\rightarrow x,
\,t\downarrow0}\frac{\phi(z+th)-\phi(z)}{t}.
$$

Recall that $\phi$ is said to be regular at $x$ if $\phi$ is local Lipschitz at $x$  and admits directional derivatives $\phi'(x,h)$ at $x$ for all $h\in X$ with $\phi'(x,h)=\phi^{\circ}(x,h)$, where $\phi'(x,h)$ is defined by
$$
\phi'(x,h):=\lim_{t\rightarrow 0^+}\frac{\phi(x+th)-\phi(x)}{t}.
$$

The following lemma, cited from \cite[Proposition 2.1.1]{C}, is useful in our analysis.

\begin{lem} Let  $\phi :X\rightarrow \mathbb{R}\cup\{+\infty\}$ be a proper
lower semicontinuous function and $x\in \mathrm{dom}(\phi)$. Suppose that $\phi$ is local
Lipschitz at $x$. Then

{\rm(i)} $\phi^{\circ}(x,\cdot)$ is a finite, positively homogeneous and subadditive function on $X$.

{\rm(ii)} $\partial_c\phi(x)$ is a nonempty weak$^*$-compact subset.

{\rm(iii)} $\phi^{\circ}(x,h)=\max\{\langle x^*,h\rangle: x^*\in\partial_c\phi(x)\}$ for all $h\in X$.

\end{lem}

Let $C$ be a convex subset of $X$. The set $[0,1]C$ is defined by
\begin{equation}\label{2.3}
[0,1]C:=\left\{
\begin{array}r
\{tx: t\in [0,1]\ {\rm and} \ x\in C\},\ \ {\rm if} \ C\not=\emptyset,\\
\{0\}\ \ \ \ \ \ \ \ \ \ \ \ \ \ \ \ , \ \ {\rm if} \ C=\emptyset.
\end{array}
\right.
\end{equation}
Recall from \cite{H} that the end set of $C$ is defined as
\begin{equation}\label{2.4}
  E[C]:=\{z\in\overline{[0, 1]C}:  tz\not\in\overline{[0, 1]C} \ {\rm for\ all} \ t>1\}.
\end{equation}
It is known from \cite[Lemma 1.1]{H1} that
\begin{equation}\label{2.4}
  E[C]=\{z\in\overline C: tz\not\in\overline C \ {\rm for\ all} \ t>1\}.
\end{equation}

We end this section with the following lemma on the end set which is cited from \cite[Lemmas 2.1 and 2.2]{H}

\begin{lem}
Let $C$ be a convex subset in $X$. Then
\begin{itemize}
    \item[\rm(i)]  $0\not\in E[C]$.

\item[\rm(ii)]  If $C$ is closed and contains the origin, then $E[C]=\{z\in C: tz\not\in C\ \forall t>1\}$ and $E[C]\subset C$.

\item[\rm(iii)]  Suppose that there exist $t_0>0$ and $z_0\in C$ such that $0\not=t_0z_0\in C$. If $M:=\sup\{t>0:tz_0\in\overline{[0,1]C}\}<+\infty$, then $M z_0\in E[C]$.

\item[\rm(iv)] $E[C]=\emptyset$ if $C$ is a cone.
\end{itemize}

\end{lem}

\setcounter{equation}{0}

\section{Constraint qualifications by Clarke normal cone and subdifferential}

In this section, we study constraint qualifications involving Clarke normal cone and subdifferential for a nonconvex inequality and target at sufficient and/or necessary conditions ensuring these qualifications.\\

\noindent{\bf 3.1 Constraint qualifications for local Lipschitz functions}\\

Throughout this subsection, we suppose that $\phi:X\rightarrow \mathbb{R}$ is a local Lipschitz function. We consider the following inequality:
\begin{equation}\label{3.1}
\phi(x)\leq 0.
\end{equation}
We denote by $S:=\{x\in X: \phi(x)\leq 0\}$ the solution set of inequality (3.1).

We discuss the following constraint qualifications for inequality (3.1) that are given by Clarke subdifferential and normal cone.

{\it Let $\bar x\in {\rm bd}(S)$ and $\tau>0$. We say that

{\rm(i)} inequality (3.1) satisfies Clarke BCQ at $\bar x$ if
\begin{equation}\label{3.2a}
  N_c(S,\bar x)\subset[0,+\infty)\partial_c\phi(\bar x);
\end{equation}

{\rm(ii)} inequality (3.1) satisfies Clarke $\tau$-strong BCQ at $\bar x$ if}
\begin{equation}\label{3.3a}
  N_c(S,\bar x)\cap B_{X^*}\subset[0,\tau]\partial_c\phi(\bar x).
\end{equation}

\noindent{\bf Remark 3.1} For the case when $\phi$ is continuous and convex, it is trivial that $\partial\phi(\bar x)\subset N(S,\bar x)$ and thus \eqref{3.2a} is equivalent to BCQ for convex inequality \eqref{1.1}. Further, by Clarke normal cone and subdifferential, it is not feasible for extending BCQ to the nonconvex inequality \eqref{3.1} via the following equation form:
\begin{equation}\label{3.4b}
 N_c(S,\bar x)=[0,+\infty)\partial_c\phi(\bar x),
\end{equation}
since the inclusion $\partial_c\phi(\bar x)\subset N_c(S,\bar x)$ does not hold trivially. For example, let $X=\mathbb{R}$, $\phi(x)=-|x|$ and $\bar x=0$. Then by computation, one has
$$N_c(S,\bar x)=N_c(X,\bar x)=\{0\}\ \ {\rm and} \ \ \partial_c\phi(\bar x)=[-1, 1],
$$
which shows that $\partial_c\phi(\bar x)\not\subset N_c(S,\bar x)$. \hfill $\Box$\\

By the concept of end set in \eqref{2.4}, the following theorem characterizes the difference between Clarke strong BCQ and Clarke BCQ. This result also shows that Clarke strong BCQ is stronger than Clarke BCQ.

\begin{them}
Let $\bar x\in {\rm bd}(S)$ and $\tau>0$. Then inequality \eqref{3.1} satisfies Clarke $\tau$-BCQ at $\bar x$ if and only if inequality \eqref{3.1} satisfies Clarke BCQ at $\bar x$ and
\begin{equation}\label{3.2}
d\big(0, E[\partial_c\phi(\bar x)\cap N_c(S,\bar x)]\big)\geq\frac{1}{\tau}.
\end{equation}

\end{them}

{\bf Proof.} The necessity part. Since Clarke strong BCQ implies Clarke BCQ trivially, it suffices to prove \eqref{3.2}. To do this, let $x^*\in E[\partial_c\phi(\bar x)\cap N_c(S,\bar x)]$. Then $x^*\not=0$ and
$$
x^*\in \overline{[0,1](\partial_c\phi(\bar x)\cap N_c(S,\bar x))}=([0,1]\partial_c\phi(\bar x))\cap N_c(S,\bar x)\subset N_c(S,\bar x)
$$
(thanks to Lemmas 2.1 and 2.2). This and Clarke $\tau$-strong BCQ imply that
$$
\frac{x^*}{\|x^*\|}\in N_c(S,\bar x)\cap B_{X^*}\subset [0,\tau]\partial_c\phi(x).
$$
Thus, there exist $t\in (0,\tau]$ and $z^*\in\partial_c\phi(x)$ such that $\frac{x^*}{\|x^*\|}=tz^*$ and consequently
$$z^*=\frac{x^*}{t\|x^*\|}\in\partial_c\phi(x)\cap N_c(S,\bar x)\subset ([0,1]\partial_c\phi(\bar x))\cap N_c(S,\bar x).
 $$
Noting that $x^*\in E[\partial_c\phi(\bar x)\cap N_c(S,\bar x)]$, it follows from the definition of the end set in \eqref{2.4} that $\frac{1}{t\|x^*\|}\leq 1$, which implies that
$$\|x^*\|\geq\frac{1}{t}\geq\frac{1}{\tau}.
$$
Hence \eqref{3.2} holds.

The sufficiency part. Let $x^*\in N_c(S,\bar x)\cap B_{X^*}$ with $\|x^*\|>0$. By Clarke BCQ of \eqref{3.2} at $\bar x$, there exist $t_0>0$ and $x_0^*\in\partial_c\phi(\bar x)$ such that $x^*=tx_0^*$. Set $$M:=\sup\{t>0:tx^*\in ([0,1]\partial_c\phi(\bar x))\cap N_c(S,\bar x)\}.$$
Then $0<M<+\infty$ by Lemma 2.1 and so it follows from Lemma 2.2 that
\begin{equation}\label{3.5}
M x^*\in E[\partial_c\phi(\bar x)\cap N_c(S,\bar x)]\subset ([0,1]\partial_c\phi(\bar x))\cap N_c(S,\bar x).
\end{equation}
By virtue of \eqref{3.2}, one has
$$
M\geq\|M x^*\|\geq d(0, E[\partial_c\phi(\bar x)\cap N_c(S,\bar x)])\geq\frac{1}{\tau}.
$$
This and \eqref{3.5} imply that
$$
x^*\in [0,\frac{1}{M}]\partial_c\phi(\bar x)\subset[0,\tau]\partial_c\phi(\bar x).
$$
Hence Clarke $\tau$-strong BCQ of \eqref{3.3a} holds at $\bar x$. The proof is complete.\hfill$\Box$\\

For the case that $\partial_c\phi(\bar x)\subset N_c(S,\bar x)$, we have the following sharper corollary on Clarke strong BCQ and Clarke BCQ of inequality \eqref{3.1}. The proof follows immediately from Theorem 3.1.
\begin{coro}
Let $\bar x\in {\rm bd}(S)$ and $\tau>0$. Suppose that $\partial_c\phi(\bar x)\subset N_c(S,\bar x)$. Then inequality \eqref{3.1} satisfies Clarke $\tau$-BCQ at $\bar x$ if and only if inequality \eqref{3.1} satisfies Clarke BCQ at $\bar x$ and
\begin{equation}\label{3.6c}
d\big(0, E[\partial_c\phi(\bar x)]\big)\geq\frac{1}{\tau}.
\end{equation}

\end{coro}

\noindent{\bf Remark 3.2} It is shown by Theorem 3.1 that Clarke strong BCQ at $\bar x\in {\rm bd}(S)$ is valid if and only if Clarke subdifferntial $\partial_c\phi(\bar x)$ contains all directions of $N_c(S,\bar x)$ and the end set of $\partial_c\phi(\bar x)\cap N_c(S,\bar x)$ is strictly separated from the origin. When restricted to the case that $\phi$ is a continuous convex function, Corollary 3.1 reduces to \cite[Theorem 3.1]{H} since $\partial\phi(\bar x)\subset N(S,\bar x)$ holds trivially. \hfill$\Box$\\

The following proposition is on the condition $\partial_c\phi(\bar x)\subset N_c(S,\bar x)$.
\begin{pro}
Let $\bar x\in {\rm bd}(S)$. Then
\begin{equation}\label{3.6a}
\partial_c\phi(\bar x)\subset N_c(S,\bar x) \Longleftrightarrow T_c(S,\bar x)\subset\{h\in X: \phi^{\circ}(\bar x,h)\leq 0\}.
\end{equation}
Suppose further that $\phi$ is regular at $\bar x$. Then $\partial_c\phi(\bar x)\subset N_c(S,\bar x)$.
\end{pro}

{\bf Proof.} By Lemma 2.1, one has that $\partial_c\phi(\bar x)$ is weak$^*$-closed and thus
\begin{equation}\label{3.8b}
 \partial_c\phi(\bar x)\subset N_c(S,\bar x) \Longleftrightarrow \sigma_{\partial_c\phi(\bar x)}\leq\sigma_{N_c(S,\bar x)},
\end{equation}
where $\sigma_U$ is the supporting function of $U$. By virtue of Lemma 2.1 again, for any $h\in X$, one has
\begin{equation*}
\sigma_{\partial_c\phi(x)}(h)=\phi^{\circ}(\bar x,h) \ \ {\rm and} \ \
\sigma_{N_c(S,\bar x)}(h)=\left\{
\begin{array}l
\ 0, \ \ \ \ h\in T_c(S,\bar x),\\
+\infty,\ h\not\in T_c(S,\bar x).
\end{array}
\right.
\end{equation*}
This and \eqref{3.8b} imply that the equivalence in \eqref{3.6a} follows.

Suppose that $\phi$ is regular at $\bar x$. We only need to show that
\begin{equation}\label{3.5b}
  T_c(S,\bar x)\subset\{h\in X: \phi^{\circ}(\bar x,h)\leq 0\}.
\end{equation}

Let $h\in T_c(S,\bar x)$. Since $\phi$ is regular at $\bar x$, it follows that there exists $t_n\rightarrow 0^+$ such that
\begin{equation}\label{3.6b}
  \phi^{\circ}(\bar x ,h)=\phi'(\bar x,h)=\lim_{n\rightarrow \infty}\frac{\phi(\bar x+t_nh)-\phi(\bar x)}{t_n}.
\end{equation}
Note that $h\in T_c(S,\bar x)$ and thus there exists $h_n\rightarrow h$ such that $\bar x+t_nh_n\in S$ for each $n$. Noting that $\phi$ is local Lipschitz around $\bar x\in {\rm bd}(S)$, it follows that $\phi(\bar x)=0$ and
\begin{equation*}
\lim_{n\rightarrow \infty}\frac{\phi(\bar x+t_nh)-\phi(\bar x)}{t_n}=\limsup_{n\rightarrow \infty}\frac{\phi(\bar x+t_nh_n)}{t_n}\leq 0.
\end{equation*}
This and \eqref{3.6b} imply that $\phi^{\circ}(\bar x ,h)\leq 0 $ and so \eqref{3.5b} holds. The proof is complete. \hfill$\Box$\\

The following corollary is immediate from Theorem 3.1 and Corollary 3.1.
\begin{coro}
Let $\bar x\in {\rm bd}(S)$ and $\tau>0$. Suppose that $\phi$ is regular at $\bar x$. Then inequality \eqref{3.1} satisfies Clarke $\tau$-BCQ at $\bar x$ if and only if inequality \eqref{3.1} satisfies Clarke BCQ at $\bar x$ and \eqref{3.6c} holds.
\end{coro}

We are now in a position to study Clarke BCQ and Clarke strong BCQ for inequality \eqref{3.1} and pay attention to several necessary and/or sufficient conditions for these constraint qualifications via the following theorems.

\begin{them}
Let $\bar x\in {\rm bd}(S)$.

{\rm(i)} Suppose that inequality \eqref{3.1} satisfies Clarke BCQ at $\bar x$. Then
\begin{equation}\label{3.11a}
  \{h\in X: \phi^{\circ}(\bar x,h)\leq 0\}\subset T_c(S,\bar x).
\end{equation}

{\rm(ii)} Suppose that $0\not\in\partial_c\phi(\bar x)$. Then inequality \eqref{3.1} satisfies Clarke BCQ at $\bar x$ if and only if  \eqref{3.11a} holds.

\end{them}

{\bf Proof.} (i) It is trivial that Clarke BCQ of \eqref{3.1} at $\bar x$ implies that
\begin{equation}\label{3.12b}
   N_c(S,\bar x)\subset {\rm cl}^*([0,+\infty)\partial_c\phi(\bar x)),
\end{equation}
where ${\rm cl}^*(\cdot)$ denotes the closure with respect to weak$^*$ topology. We next prove that \eqref{3.11a}$\Leftrightarrow$\eqref{3.12b}.

It is known that \eqref{3.12b} is equivalent to
\begin{equation}\label{3.6}
\sigma_{N_c(S,\bar x)}\leq\sigma_{{\rm cl}^*([0,+\infty)\partial_c\phi(\bar x))}=\sigma_{[0,+\infty)\partial_c\phi(\bar x)}.
\end{equation}
By computation, one has
\begin{equation*}
\sigma_{N_c(S,\bar x)}=\left\{
\begin{array}l
\ 0, \ \ \ \ h\in T_c(S,\bar x),\\
+\infty,\ h\not\in T_c(S,\bar x),
\end{array}
\right.
{\rm and} \ \
\sigma_{[0,+\infty)\partial_c\phi(x))}=\left\{
\begin{array}l
\ 0,\ \ \ \ {\rm if} \  \phi^{\circ}(\bar x,h)\leq 0,\\
+\infty,\ \, {\rm if} \ \phi^{\circ}(\bar x,h)> 0.
\end{array}
\right.
\end{equation*}
This means that \eqref{3.6} holds if and only if
$$
\{h\in X: \phi^{\circ}(\bar x,h)\leq 0\}\subset T_c(S,\bar x).
$$

(ii) To show the equivalence between Clarke BCQ  and \eqref{3.11a}, we only need to prove that $[0,+\infty)\partial_c\phi(x)$ is weak$^*$-closed by (i).

To do this, let $x^*\in{\rm cl}^*([0,+\infty)\partial_c\phi(\bar x))$. Then there exist nets $t_\alpha\geq 0$ and $x^*_\alpha\in \partial_c\phi(\bar x)$ such that $t_\alpha x_\alpha^*\stackrel{w^*}\rightarrow x^*$. We claim that $\{t_\alpha^*\}$ is bounded. (Otherwise, we can assume that $t_\alpha\rightarrow +\infty$ and thus $x_\alpha^*=\frac{1}{t_\alpha}t_\alpha x_\alpha^*\stackrel{w^*}\rightarrow 0$. This means that $0\in\partial_c\phi(\bar x)$ by Lemma 2.1, which contradicts $0\not\in\partial_c\phi(\bar x)$). By virtue of Lemma 2.1 again, without loss of generality, we can assume that $x^*_\alpha\stackrel{w^*}\rightarrow x_0^*\in \partial_c\phi(\bar x)$ and $t_\alpha^*\rightarrow t_0\geq 0$ (considering subnet if necessary). Then
$$
x^*=t_0x_0^*\in [0,+\infty)\partial_c\phi(\bar x).
$$
Hence $[0,+\infty)\partial_c\phi(\bar x)$ is weak$^*$-closed.  The proof is complete.\hfill$\Box$\\

The following theorem provides a characterization for Clarke strong BCQ of inequality \eqref{3.1}.

\begin{them}
Let $\bar x\in {\rm bd}(S)$ and $\tau>0$.  Then inequality \eqref{3.1} satisfies Clarke $\tau$-strong BCQ at $\bar x$ if and only if
\begin{equation}\label{3.7}
d(h,T_c(S,\bar x))\leq \tau\max\{0,\phi^{\circ}(\bar x,h)\}
\end{equation}
holds for any $h\in X$.
\end{them}

{\bf Proof.} It is known that Clarke $\tau$-strong BCQ of \eqref{3.3a} at $\bar x$ is equivalent to
\begin{equation}\label{3.8}
\sigma_{N_c(S,\bar x)\cap B_{X^*}}\leq\sigma_{[0,\tau]\partial_c\phi(\bar x)}.
\end{equation}
Noting that $\sigma_{N_c(S,\bar x)\cap B_{X^*}}=d(\cdot, T_c(S,\bar x))$ and by Lemma 2.1, one has
$$
\sigma_{[0,\tau]\partial_c\phi(x)(h)}=\tau\max\{0,\phi^{\circ}(\bar x,h)\}\ \ \forall h\in X.
$$
This and \eqref{3.8} imply that Clarke $\tau$-strong BCQ  of \eqref{3.3a} at $\bar x$ holds if and only if \eqref{3.7} holds for any $h\in X$. The proof is complete.\hfill$\Box$\\

\noindent{\bf Remark 3.3} When $\phi$ is a convex continuous function, it is known that Clarke directional derivative $\phi^{\circ}(\bar x,h)$ coincide with directional derivative $\phi'(\bar x,h)$ and then Theorem 3.3 reduces to that the validity of $\tau$-strong BCQ at $\bar x$ is equivalent to
$$d(h, T(S,\bar x))\leq \tau\max\{0,\phi'(\bar x,h)\}\ \  \forall h\in X.
$$
This result is also proved by Zheng and Ng \cite{Z1}. Readers can refer to \cite[Theorem 2.3]{Z1} for more details and its proof.\hfill$\Box$\\

Suppose that $\varphi:X\rightarrow \mathbb{R}\cup\{+\infty\}$ is a proper convex lower semicontinuous function and $\tau\in (0, +\infty)$. Recall that convex inequality $\varphi(x)\leq 0$ is said to have the global error bound with constant $\tau>0$ if
\begin{equation}\label{3.14}
  d(x,S_{\varphi})\leq\tau\max\{0,\varphi(x)\}\ \ \forall x\in X,
\end{equation}
where $S_{\varphi}:=\{u\in X: \varphi(u)\leq 0\}$.\\

The following theorem provides conditions for ensuring Clarke strong BCQ of inequality \eqref{3.1}. The proof mainly relies on Theorem 3.3.

\begin{them}
Let $\bar x\in {\rm bd}(S)$ and $\tau>0$.

{\rm (i)} Suppose that inequality \eqref{3.1} satisfies Clarke BCQ at $\bar x$ and convex inequality $\phi^{\circ}(\bar x,h)\leq 0$ has the global error bound with constant $\tau>0$. Then inequality \eqref{3.1} satisfies Clarke $\tau$-strong BCQ at $\bar x$.

{\rm (ii)} Suppose that $\partial_c\phi(\bar x)\subset N_c(S,\bar x)$. Then inequality \eqref{3.1} satisfies Clarke $\tau$-strong BCQ at $\bar x$ if and only if inequality \eqref{3.1} satisfies Clarke BCQ at $\bar x$ and convex inequality $\phi^{\circ}(\bar x,h)\leq 0$ has the global error bound with constant $\tau>0$.
\end{them}

{\bf Proof.} We denote
$$
S_{\phi^{\circ}(\bar x,\cdot)}:=\{h\in X: \phi^{\circ}(\bar x,h)\leq 0\}.
$$

(i) Since convex inequality $\phi^{\circ}(\bar x,h)\leq 0$ has the global error bound with constant $\tau>0$, it follows that
\begin{equation}\label{3.15}
d(h,S_{\phi^{\circ}(\bar x,\cdot)})\leq \tau\max\{0,\phi^{\circ}(\bar x,h)\} \  \ \forall h\in X.
\end{equation}
Noting that inequality \eqref{3.1} satisfies Clarke BCQ at $\bar x$, it follows from Theorem 3.2 that \eqref{3.11a} holds; that is, $S_{\phi^{\circ}(\bar x,\cdot)}\subset T_c(S,\bar x)$. This and \eqref{3.15} imply that
\begin{equation*}
d(h,T_c(S,\bar x))\leq \tau\max\{0,\phi^{\circ}(\bar x,h)\} \ \ \forall h\in X.
\end{equation*}
Hence Clarke $\tau$-strong BCQ of inequality \eqref{3.1} at $\bar x$ follows from Theorem 3.3.

(ii) Since Clarke BCQ of \eqref{3.1} at $\bar x$ follows from Clarke $\tau$-strong BCQ trivially, we only need to prove the global error bound for convex inequality $\phi^{\circ}(\bar x,h)\leq 0$ holds with constant $\tau>0$.

By virtue of Theorem 3.3, one has
\begin{equation}\label{3.14c}
d(h,T_c(S,\bar x))\leq \tau\max\{0,\phi^{\circ}(\bar x,h)\} \  \ \forall h\in X.
\end{equation}
Note that $\partial_c\phi(\bar x)\subset N_c(S,\bar x)$ and so Proposition 3.1 implies that $T_c(S,\bar x)\subset S_{\phi^{\circ}(\bar x,\cdot)}$.  This and \eqref{3.14c} imply that
\begin{equation*}
d(h,S_{\phi^{\circ}(\bar x,\cdot)})\leq \tau\max\{0,\phi^{\circ}(\bar x,h)\} \  \ \forall h\in X.
\end{equation*}
Hence $\phi^{\circ}(\bar x,h)\leq 0$ has global error bound with $\tau>0$. The proof is complete.\hfill$\Box$\\

The following corollary is immediate from Proposition 3.1 and Theorem 3.4.

\begin{coro}
Let $\bar x\in {\rm bd}(S)$ and $\tau>0$. Suppose that $\phi$ is regular at $\bar x$. Then inequality \eqref{3.1} satisfies Clarke $\tau$-strong BCQ at $\bar x$ if and only if inequality \eqref{3.1} satisfies Clarke BCQ at $\bar x$ and convex inequality $\phi^{\circ}(\bar x,h)\leq 0$ has the global error bound with constant $\tau>0$.
\end{coro}

\noindent{\bf 3.2 Constraint qualifications  for lower semicontinuous functions}\\

This subsection is devoted to the  study of constraint qualifications for the nonconvex inequality which is defined by a proper lower semicontinuous function. For this aim, we suppose that $\phi:X\rightarrow \mathbb{R}\cup\{+\infty\}$ is a proper lower semicontinuous function throughout this subsection.

For the case that $\phi$ is convex, Zheng and Ng \cite{Z1} studied the extend BCQ and strong BCQ for convex inequality \eqref{3.1}.

{\it Let $\bar x\in {\rm bd}(S)\cap \phi^{-1}(0)$. Recall from \cite{Z1} that convex inequality \eqref{3.1} is said to satisfy the extended BCQ at $\bar x$ if
\begin{equation}\label{4.2a}
  N(S,\bar x)=[0,+\infty)\partial\phi(\bar x)+\partial^{\infty}\phi(\bar x)
\end{equation}
and convex inequality \eqref{3.1} is said to satisfy the strong BCQ at $\bar x$ if there exists $\tau>0$ such that}
\begin{equation}\label{4.3a}
N(S,\bar x)\cap B_{X^*}\subset [0, \tau]\partial\phi(\bar x)+\partial^{\infty}\phi(\bar x).
\end{equation}

When $\phi$ is a continuous convex function, it follows that $\partial^{\infty}\phi(\bar x)=\{0\}$ and thus the extended BCQ and strong BCQ reduce to those in \eqref{1.1} and \eqref{1.2}, respectively.

The main work in this subsection is to consider two types of constraint qualifications for inequality \eqref{3.1} with the help of Clarke singular subdifferential. This is inspired by basic constraint qualifications in \eqref{4.2a} and \eqref{4.3a}.

We first present following proposition related to some useful properties of Clarke subdifferential and Clarke singular subdifferential. This result is an extension of \cite[Lemma 4.1]{H} to the nonconvex case.

\begin{pro}
Let $x\in{\it dom}(\phi)$ be such that $\partial_c\phi(x)\not=\emptyset$ and $r>0$. Then

{\rm(i)} $\partial_c\phi(x)+r\partial_c^{\infty}\phi(x)=\partial_c\phi(x)$.

{\rm(ii)}  $(0,r]\partial_c\phi(x)+\partial_c^{\infty}\phi(x)=(0,r]\partial_c\phi(x)$.

{\rm(iii)}  $\partial_c^{\infty}\phi(x)\subset\overline{(0,r]\partial_c\phi(x)}$.

{\rm(iv)}  $\overline{[0,r]\partial_c\phi(x)}=[0,r]\partial_c\phi(x) +\partial_c^{\infty}\phi(x)$.

\end{pro}

{\bf Proof.} (i) Note that $0\in\partial_c^{\infty}\phi(x)$ and thus the inverse inclusion holds. Conversely, let $x^*\in \partial_c\phi(x)$ and $u^*\in\partial_c^{\infty}\phi(x)$. Then
$$
(x^*,-1)\in N_c({\rm epi}(\phi), (x,\phi(x)))\ \  {\rm and} \ \ (u^*,0)\in N_c({\rm epi}(\phi), (x,\phi(x))).
$$
Since $N_c({\rm epi}(\phi), (x,\phi(x)))$ is a closed convex cone, it follows that
\begin{eqnarray*}
(x^*+ru^*,-1)&\in& N_c({\rm epi}(\phi), (x,\phi(x)))+N_c({\rm epi}(\phi), (x,\phi(x)))\\
&=&N_c({\rm epi}(\phi), (x,\phi(x))).
\end{eqnarray*}
This implies that $x^*+ru^*\in \partial_c\phi(x)$.

(ii) Since $0\in\partial_c^{\infty}\phi(x)$, we have
$$(0,r](\partial_c\phi(x))+\partial_c^{\infty}\phi(x)\supset(0,r]\partial_c\phi(x).
$$
Conversely, for any $(t,x^*,u^*)\in (0,r]\times \partial_c\phi(x)\times \partial_c^{\infty}\phi(x)$, by (i), one has $x^*+\frac{1}{t}u^*\in \partial_c\phi(x)$ and thus
$$tx^*+u^*=t(x^*+\frac{1}{t}u^*)\in (0,r]\partial_c\phi(x).
$$

(iii) Let $u^*\in\partial_c^{\infty}\phi(x)$. Take any $x^*\in\partial_c\phi(x)$ and set $u_n^*:=\frac{1}{n}x^*+u^*$ for each $n\in\mathbb{N}$. Then for any $n$ sufficiently large, one has
\begin{eqnarray*}
u_n^*=\frac{1}{n}x^*+u^*&\in&(0,r]\partial_c\phi(x)+\partial_c^{\infty}\phi(x) \\ &=&(0,r](\partial_c\phi(x))
\end{eqnarray*}
(the equation follows from (ii)). Thus $u^*\in\overline{(0,r]\partial_c\phi(x)}$.

(iv) Let $z^*\in \overline{[0,r](\partial_c\phi(x))}$. Then there exist $t_n\in [0,r]$ and $x_n^*\in\partial_c\phi(x)$ such that $t_nx_n^*\rightarrow z^*$. Without loss of generalization, we can assume that $t_n\rightarrow t_0\in[0,r]$ (considering subsequence if necessary). Note that
$$
(x_n^*,-1)\in N_c({\rm epi}(\phi), (x,\phi(x)))
$$
and so $(z^*,-t_0)\in N_c({\rm epi}(\phi), (x,\phi(x)))$. We divided $t_0$ into two cases:

Case 1: $t_0=0$. Then $(z^*,0)\in N_c({\rm epi}(\phi), (x,\phi(x)))$ and thus
$$
z^*\in \partial_c^{\infty}\phi(x)\subset [0,r](\partial_c\phi(x))+\partial_c^{\infty}\phi(x).
$$

Case 2: $t_0>0$. Then $t_0\in(0,r]$ and thus
$$
(\frac{1}{t_0}z^*,-1)\in N_c({\rm epi}(\phi), (x,\phi(x))).
$$
This implies that $\frac{1}{t_0}z^*\in\partial_c\phi(x)$ and consequently
\begin{eqnarray*}
z^*=t_0(\frac{1}{t_0}z^*)&\in& (0,r](\partial_c\phi(x))\\
&\subset& [0,r](\partial_c\phi(x))+\partial_c^{\infty}\phi(x).
\end{eqnarray*}

Conversely, let $y^*:=tx^*+u^*\in[0,r](\partial_c\phi(x))+\partial_c^{\infty}\phi(x)$ where $t\in [0,r]$ , $x^*\in \partial_c\phi(x)$ and $u^*\in\partial_c^{\infty}\phi(x)$. If $t=0$, then $y^*=u^*\in \partial_c^{\infty}\phi(x)\subset \overline{[0,r](\partial_c\phi(x))}$ by (iii). If $t>0$, then (ii) implies that
\begin{eqnarray*}
y^*&\in& (0,r](\partial_c\phi(x))+\partial_c^{\infty}\phi(x)\\
&=&(0,r](\partial_c\phi(x))\\
&\subset&\overline{[0,r](\partial_c\phi(x))}.
\end{eqnarray*}
The proof is complete. \hfill$\Box$\\

By Proposition 3.2, the main result in this subsection is given via the following theorem.

\begin{them}
Let $\bar x\in{\rm bd}(S)\cap\phi^{-1}(0)$ and $\tau>0$. Then the following statements are equivalent:
\begin{itemize}
\item[{\rm(i)}]  The extended Clarke $\tau$-strong BCQ at $\bar x$ holds:
\begin{equation}\label{4.2}
 N_c(S,\bar x)\cap B_{X^*}\subset [0,\tau]\partial_c\phi(\bar x)+\partial_c^{\infty}\phi(\bar x);
\end{equation}

\item[{\rm(ii)}]  The extended Clarke  BCQ at $\bar x$ holds:
\begin{equation}\label{4.3}
  N_c(S,\bar x)\subset [0,+\infty)\partial_c\phi(\bar x)+\partial_c^{\infty}\phi(\bar x)
\end{equation}
and
\begin{equation}\label{4.4}
  d\big(0, E\big[([0,1]\partial_c\phi(\bar x)+\partial_c^{\infty}\phi(\bar x))\cap N_c(S,\bar x)\big]\big)\geq\frac{1}{\tau}.
\end{equation}
\end{itemize}
\end{them}

{\bf Proof.} (i)$\Rightarrow$(ii): We divided $\partial_c\phi(\bar x)$ into two cases:

Case 1: $\partial_c\phi(\bar x)=\emptyset$. Then $N_c(S,\bar x)\cap B_{X^*}\subset \partial_c^{\infty}\phi(\bar x)$ and thus \eqref{4.3} holds. Noting that $(\partial_c^{\infty}\phi(\bar x)\cap N_c(S,\bar x)$ is a cone, it follows from Lemma 2.2 that
$$ E[([0,1]\partial_c\phi(\bar x)+\partial_c^{\infty}\phi(\bar x))\cap N_c(S,\bar x)]=\emptyset.
$$
This implies that
\begin{equation*}
  d\big(0, E\big[([0,1]\partial_c\phi(\bar x)+\partial_c^{\infty}\phi(\bar x))\cap N_c(S,\bar x)\big]\big)=+\infty.
\end{equation*}
Thus \eqref{4.4} hold.

Case 2:  $\partial_c\phi(\bar x)\not=\emptyset$. Note that $\eqref{4.3}$ follows from \eqref{4.2} and we next prove \eqref{4.4}. To do this, let $x^*\in E[([0,1]\partial_c\phi(\bar x)+\partial_c^{\infty}\phi(\bar x))\cap N_c(S,\bar x)]$. By Proposition 3.2 and Lemma 2.2, one has
$$
0\not=x^*\in ([0,1]\partial_c\phi(\bar x)+\partial_c^{\infty}\phi(\bar x))\cap N_c(S,\bar x)\subset N_c(S,\bar x).
$$
Then \eqref{4.2} implies that
$$
\frac{x^*}{\|x^*\|}\in N_c(S,\bar x)\cap B_{X^*}\subset [0,\tau]\partial_c\phi(\bar x)+\partial_c^{\infty}\phi(\bar x).
$$
Thus, there exist $t\in [0,\tau]$, $y^*\in \partial_c\phi(\bar x)$ and $u^*\in\partial_c^{\infty}\phi(\bar x)$ such that
\begin{equation}\label{4.5}
\frac{x^*}{\|x^*\|}=ty^*+u^*.
\end{equation}

We claim that $t>0$. Suppose on the contrary that $t=0$. Then \eqref{4.5} implies that
$$
x^*=\|x^*\|u^*\in \partial_c^{\infty}\phi(\bar x)
$$
Noting that $\partial_c^{\infty}\phi(\bar x)$ is a cone, it follows from Proposition 3.2(iii) that
$$
\lambda x^*\in\partial_c^{\infty}\phi(\bar x)\subset [0,1]\partial_c\phi(\bar x)+\partial_c^{\infty}\phi(\bar x) \ \  \forall \lambda>0,
$$
which contradicts $x^*\in E[([0,1]\partial_c\phi(\bar x)+\partial_c^{\infty}\phi(\bar x))\cap N_c(S,\bar x)]$.

Since $y^*+\frac{1}{t}u^*\in \partial_c\phi(\bar x)+\frac{1}{t}\partial_c^{\infty}\phi(\bar x)$, it follows that
$$
\frac{x^*}{t\|x^*\|}=y^*+\frac{1}{t}u^*\in ([0,1]\partial_c\phi(\bar x)+\partial_c^{\infty}\phi(\bar x))\cap N_c(S,\bar x)
$$
and consequently $\frac{1}{t\|x^*\|}\leq 1$ by $x^*\in E[([0,1]\partial_c\phi(\bar x)+\partial_c^{\infty}\phi(\bar x))\cap N_c(S,\bar x)]$. This implies that $$\|x^*\|\geq\frac{1}{t}\geq\frac{1}{\tau}.
$$
Hence \eqref{4.4} holds.

(i)$\Rightarrow$(ii): If $\partial_c\phi(\bar x)=\emptyset$ or $\partial_c\phi(\bar x)=\{0\}$, then \eqref{4.3} reduces to $N_c(S,\bar x)\subset \partial_c\phi(\bar x)$ by \eqref{2.3}  and thus
$$
N_c(S,\bar x)\cap B_{X^*}\subset \partial_c^{\infty}\phi(\bar x)=[0,\tau]\partial_c\phi(\bar x)+\partial_c^{\infty}\phi(\bar x).
$$

We next suppose that there exists $z^*\in\partial_c\phi(\bar x)\backslash\{0\}$. Let $x^*\in N_c(S,\bar x)\cap B_{X^*}$ be such that $\|x^*\|>0$. By \eqref{4.3}, there exist $t\geq 0$, $u^*\in\partial_c\phi(\bar x)$ and $v^*\in\partial_c\phi^{\infty}(\bar x)$ such that
\begin{equation}\label{4.6}
  x^*=tu^*+v^*.
\end{equation}
If $t=0$, then $x^*=v^*\in \partial_c^{\infty}\phi(\bar x)=[0,\tau]\partial_c\phi(\bar x)+\partial_c^{\infty}\phi(\bar x)$ and the proof is completed.

We next consider $t>0$. By \eqref{4.6} and Proposition 3.2, one has
$$
\frac{x^*}{t}=u^*+\frac{v^*}{t}\in\partial_c\phi(\bar x)+\frac{1}{t}\partial_c^{\infty}\phi(\bar x)=\partial_c\phi(\bar x)+\partial_c^{\infty}\phi(\bar x)\subset [0,1]\partial_c\phi(\bar x)+\partial_c^{\infty}\phi(\bar x).
$$
Let $M:=\sup\{\lambda>0: \lambda x^*\in[0,1]\partial_c\phi(\bar x)+\partial_c^{\infty}\phi(\bar x)\}$. Then $M>0$.

If $M<+\infty$, then by Lemma 2.2 one has
\begin{equation}\label{4.7}
  Mx^*\in E[([0,1]\partial_c\phi(\bar x)+\partial_c^{\infty}\phi(\bar x))\cap N_c(S,\bar x)]\subset [0,1]\partial_c\phi(\bar x)+\partial_c^{\infty}\phi(\bar x),
\end{equation}
and it follows from \eqref{4.4} that
$$
M\geq \|Mx^*\|\geq d(0, E[([0,1]\partial_c\phi(\bar x)+\partial_c^{\infty}\phi(\bar x))\cap N_c(S,\bar x)])\geq\frac{1}{\tau}.
$$
This and \eqref{4.7} imply that
$$
x^*\in [0,\frac{1}{M}]\partial_c\phi(\bar x)+\partial_c^{\infty}\phi(\bar x)\subset[0,\tau]\partial_c\phi(\bar x)+\partial_c^{\infty}\phi(\bar x).
$$

If $M=+\infty$. Then there exists $\lambda>\frac{1}{\tau}$ such that
$$
\lambda x^*\in [0,1]\partial_c\phi(\bar x)+\partial_c^{\infty}\phi(\bar x).
$$
This means that
$$
x^*\in [0,\frac{1}{\lambda}]\partial_c\phi(\bar x)+\partial_c^{\infty}\phi(\bar x)\subset[0,\tau]\partial_c\phi(\bar x)+\partial_c^{\infty}\phi(\bar x).
$$
Hence \eqref{4.2} holds.  The proof is complete.\hfill$\Box$\\

For the special case that $\partial_c\phi(\bar x)\subset N_c(S,\bar x)$, it is easy to verify that
$$
E\big[([0,1]\partial_c\phi(\bar x)+\partial_c^{\infty}\phi(\bar x))\cap N_c(S,\bar x)\big]=E[\partial_c\phi(\bar x)].
$$
Therefore the following corollary is immediate from Theorem 3.5. This result reduces to \cite[Theorem 4.1]{H} when restricted to the case that $\phi$ is a convex function

\begin{coro}
Let $\bar x\in{\rm bd}(S)\cap\phi^{-1}(0)$ and $\tau>0$. Suppose that $\partial_c\phi(\bar x)\subset N_c(S,\bar x)$. Then \eqref{4.2} holds if and only if \eqref{4.3} and $d(0,E[\partial_c\phi(\bar x)])\geq\frac{1}{\tau}$.
\end{coro}

For the case that $\partial_c^{\infty}\phi(\bar x)=\{0\}$, we have the following sharper theorem. This result reduces to Theorem 3.1 when $\phi$ is a local Lipschitz function.

\begin{them}
Let $\bar x\in{\rm bd}(S)\cap\phi^{-1}(0)$ and $\tau>0$. Suppose that $\partial_c^{\infty}\phi(\bar x)=\{0\}$. Then
\begin{equation*}
 N_c(S,\bar x)\cap B_{X^*}\subset [0,\tau]\partial_c\phi(\bar x)
\end{equation*}
holds if and only if
\begin{equation*}
  N_c(S,\bar x)\subset [0,+\infty)\partial_c\phi(\bar x)\ \ {\it and} \ \  d(0, E[\partial_c\phi(\bar x)\cap N_c(S,\bar x)])\geq\frac{1}{\tau}.
\end{equation*}
\end{them}

{\bf Proof.} Since $\partial_c^{\infty}\phi(\bar x)=\{0\}$, it follows from Proposition 4.1 that
\begin{equation}\label{4.8}
  [0,1]\partial_c\phi(\bar x)=\overline{[0,1]\partial_c\phi(\bar x)}.
\end{equation}
Note that
$$
([0,1]\partial_c\phi(\bar x))\cap N_c(S, \bar x)=[0,1](\partial_c\phi(\bar x)\cap N_c(S, \bar x))
$$
and thus
$$
[0,1](\partial_c\phi(\bar x)\cap N_c(S, \bar x))=\overline{[0,1](\partial_c\phi(\bar x)\cap N_c(S, \bar x))}
$$
(thanks to \eqref{4.8}). This and the definition of the end set in \eqref{2.4} imply that
\begin{eqnarray*}
E[\partial_c\phi(\bar x)\cap N_c(S,\bar x)]=E[([0,1]\partial_c\phi(\bar x))\cap N_c(S,\bar x)].
\end{eqnarray*}
Hence Theorem 3.6 follows from Theorem 3.5. The proof is complete.\hfill$\Box$\\

The following corollary is immediate from Corollary 3.4 and Theorem 3.6.

\begin{coro}
Let $\bar x\in{\rm bd}(S)\cap\phi^{-1}(0)$ and $\tau>0$. Suppose that $\partial_c^{\infty}\phi(\bar x)=\{0\}$ and $\partial_c\phi(\bar x)\subset N_c(S,\bar x)$. Then
\begin{equation*}
 N_c(S,\bar x)\cap B_{X^*}\subset [0,\tau]\partial_c\phi(\bar x)
\end{equation*}
holds if and only if
\begin{equation*}
  N_c(S,\bar x)\subset [0,+\infty)\partial_c\phi(\bar x)\ \ {\it and} \ \  d(0, E[\partial_c\phi(\bar x)])\geq\frac{1}{\tau}.
\end{equation*}
\end{coro}

\setcounter{equation}{0}

\section{Constraint qualifications by Fr\'echet normal cone and subdifferential}

Throughout this section, we suppose that $\phi:X\rightarrow \mathbb{R}$ is a proper lower semicontinuous function. We consider the following inequality:

\begin{equation}\label{4-1}
\phi(x)\leq 0.
\end{equation}
We still denote by $S$ the solution set of inequality \eqref{4-1}.

The main work of this section is to study constraint qualifications involving Fr\'echet subdifferential and normal cone for inequality \eqref{4-1}. For this aim, we consider the following two constraint qualifications.

{\it Let $\bar x\in {\rm bd}(S)$ and $\tau>0$. We say that

{\rm(i)} inequality \eqref{4-1} satisfies Fr\'echet BCQ at $\bar x$ if
\begin{equation}\label{4-2}
  \hat N(S,\bar x)=[0,+\infty)\hat\partial\phi(\bar x);
\end{equation}

{\rm(ii)} inequality \eqref{4-1} satisfies Fr\'echet $\tau$-strong BCQ at $\bar x$ if}
\begin{equation}\label{4-3}
  \hat N(S,\bar x)\cap B_{X^*}\subset[0,\tau]\hat\partial\phi(\bar x).
\end{equation}

It is easy to verify that constraint qualifications in \eqref{4-2} and \eqref{4-3} reduce to BCQ and strong BCQ in \cite{H,Z1} respectively when restricted to the convex case. The following proposition is useful to prove the main result of this section.

\begin{pro}
Let $\bar x\in {\rm bd}(S)$. Suppose that $\hat\partial \phi(\bar x)$ is bounded. Then
\begin{equation}\label{4-4a}
  \overline{[0, r]\hat\partial\phi(\bar x)}=[0, r]\hat\partial\phi(\bar x), \ \ \forall r>0.
\end{equation}
\end{pro}

{\bf Proof.} Let $r>0$ and $x^*\in \overline{[0, r]\hat\partial\phi(\bar x)}$ with $x^*\not=0$. Then there exist $t_n\in [0, r]$ and $x_n^*\in \hat\partial\phi(\bar x)$ such that $t_nx_n^*\rightarrow x^*$. Without loss of generalization, we can assume that $t_n\rightarrow t_0\in [0,r]$ (considering subsequence if necessary). Noting that $\hat\partial \phi(\bar x)$ is bounded and $x^*\not=0$, it follows that $t_0>0$ and thus $x_n\rightarrow \frac{1}{t_0}x^*$. This means that $\frac{1}{t_0}x^*\in \hat\partial\phi(\bar x)$ by the norm-closed property of $\hat\partial\phi(\bar x)$ and consequenlty
$$
x^*=t_0 \big(\frac{1}{t_0}x^*\big)\in [0, r]\hat\partial\phi(\bar x).
$$
Hence \eqref{4-4a} holds. The proof is complete.\hfill$\Box$\\

The following theorem provides an characterization for the Fr\'echet strong BCQ which is given by Fr\'echet BCQ and the end set of Fr\'echet subdifferential.

\begin{them}
Let $\bar x\in {\rm bd}(S)\cap \phi^{-1}(0)$ and $\tau>0$. Suppose that $\hat\partial\phi(\bar x)$ is bounded. Then inequality \eqref{4-1} satisfies Fr\'echet $\tau$-strong BCQ at $\bar x$ if and only if inequality \eqref{4-1} satisfies Fr\'echet BCQ at $\bar x$ and
\begin{equation}\label{4-5}
d\big(0, E[\hat\partial\phi(\bar x)]\big)\geq\frac{1}{\tau}.
\end{equation}

\end{them}

{\bf Proof.}
The necessity part. Note that $\hat\partial\phi(\bar x)\subset\hat N(S,\bar x)$ by the definition and thus \eqref{4-2} follows from \eqref{4-3}. We next prove \eqref{4-5}. To do this, let $x^*\in E[\hat\partial\phi(\bar x)]$. Then $x^*\not=0$ by Lemma 2.2 and thus
$$
x^*\in \overline{[0,1]\hat\partial\phi(\bar x)}=[0,1]\hat\partial\phi(\bar x)\subset \hat N(S,\bar x)
$$
(the equation follows from \eqref{4-4a}). This and Fr\'echet $\tau$-strong BCQ imply that
$$
\frac{x^*}{\|x^*\|}\in\hat N(S,\bar x)\cap B_{X^*}\subset [0,\tau]\hat\partial\phi(x).
$$
Thus, there exist $t\in (0,\tau]$ and $z^*\in\hat\partial\phi(x)$ such that $\frac{x^*}{\|x^*\|}=tz^*$ and consequently
$$
z^*=\frac{x^*}{t\|x^*\|}\in\hat\partial\phi(x)\subset [0,1]\hat\partial\phi(\bar x).
$$
Note that $x^*\in E[\hat\partial\phi(\bar x)]$ and thus $\frac{1}{t\|x^*\|}\leq 1$ by \eqref{2.4}, which implies that
$$\|x^*\|\geq\frac{1}{t}\geq\frac{1}{\tau}.
$$
Hence \eqref{4-2} holds.

The sufficiency part. Let $x^*\in \hat N(S,\bar x)\cap B_{X^*}$ with $\|x^*\|>0$. By \eqref{4-2}, there exists $t_0>0$ and $x_0^*\in\hat\partial\phi(\bar x)$ such that $x^*=tx_0^*$. We denote
$$
M:=\sup\{t>0:tx^*\in [0,1]\hat\partial\phi(\bar x)\}.
$$
Note that $\hat\partial\phi(\bar x)$ is bounded and thus $0<M<+\infty$. By Lemma 2.2, one has
\begin{equation}\label{4-7}
M x^*\in E[\hat\partial\phi(\bar x)]\subset [0,1]\hat\partial\phi(\bar x).
\end{equation}
By virtue of \eqref{4-5}, one has
$$
M\geq\|M x^*\|\geq d(0, E[\hat\partial\phi(\bar x)])\geq\frac{1}{\tau}.
$$
This and \eqref{4-7} imply that
$$
x^*\in [0,\frac{1}{M}]\partial_c\phi(\bar x)\subset[0,\tau]\hat\partial\phi(\bar x).
$$
Hence Fr\'echet $\tau$-strong BCQ of \eqref{4-3} holds at $\bar x$. The proof is complete.\hfill$\Box$\\

The following corollary is immediate from Theorem 4.1.
\begin{coro}
Let $\bar x\in {\rm bd}(S)$ and $\tau>0$. Suppose that $\phi$ is local Lipschitz at $\bar x$. Then inequality \eqref{4-1} satisfies the Fr\'echet $\tau$-strong BCQ at $\bar x$ if and only if inequality \eqref{4-1} satisfies the Fr\'echet BCQ at $\bar x$ and \eqref{4-5} holds.
\end{coro}

For the case that $X$ is a  finite dimensional space, we obtain necessary and/or sufficient conditions for Fr\'echet BCQ and strong BCQ of inequality \eqref{4-1}.

\begin{pro}
Let $\bar x\in {\rm bd}(S)$ and $\tau>0$. Suppose that $X$ is finite dimensional and $\hat\partial\phi(\bar x)$ is bounded.

{\rm(i)} If inequality \eqref{4-1} satisfies the Fr\'echet BCQ at $\bar x$, then
\begin{equation}\label{4-8}
  \big\{h\in X: \langle x^*, h\rangle\leq 0 \ \forall x^*\in\hat\partial\phi(\bar x)\big\}= \overline{\rm co}(T(S,\bar x)),
\end{equation}
where $\overline{\rm co}(U)$ is the closed convex hull of $U\subset X$.

{\rm(ii)}  Suppose that $0\not\in\hat\partial\phi(\bar x)$. Then inequality \eqref{4-1} satisfies the Fr\'echet BCQ at $\bar x$ if and only if \eqref{4-8} holds.

\end{pro}

{\bf Proof.} (i) It is known that the Fr\'echet BCQ of \eqref{4-1} at $\bar x$ implies that
\begin{equation}\label{4.8b}
  \hat N(S,\bar x)=\overline{[0,+\infty)\hat\partial\phi(\bar x)}.
\end{equation}
Note that \eqref{4.8b} is equivalent to
\begin{equation}\label{4-10}
  \sigma_{\hat N(S,\bar x)}=\sigma_{\overline{[0,+\infty)\hat \partial\phi(\bar x)}}=\sigma_{[0,+\infty)\hat \partial\phi(\bar x)}.
\end{equation}
By \eqref{2.1}, one can verify that
\begin{equation*}
  \sigma_{\hat N(S,\bar x)}(h)=\left\{
  \begin{array}r
  0,\ \ h\in \overline{\rm co}(T(S,\bar x),\\
  +\infty, h\not\in \overline{\rm co}(T(S,\bar x),
  \end{array}
  \right.
\end{equation*}
and
\begin{equation*}
  \sigma_{[0,+\infty)\hat \partial\phi(\bar x)}(h)=\left\{
  \begin{array}r
  0,\ \ {\rm if}\ \sup_{x^*\in\hat\partial\phi(\bar x)}\langle x^*,h\rangle\leq 0,\\
  +\infty, {\rm if}\ \sup_{x^*\in\hat\partial\phi(\bar x)}\langle x^*,h\rangle> 0.
  \end{array}
  \right.
\end{equation*}
This implies that \eqref{4-10} is equivalent to
$$
  \left\{h\in X: \sup_{x^*\in\hat\partial\phi(\bar x)}\langle x^*,h\rangle\leq 0 \right\}=\overline{\rm co}(T(S,\bar x)).
$$
Hence \eqref{4-8} holds.

(ii) To prove the equivalence between the Fr\'echet BCQ and \eqref{4-8}, we only need to show that $[0,+\infty)\hat \partial\phi(\bar x)$ is closed. Let $x^*\in \overline{[0,+\infty)\hat \partial\phi(\bar x)}$. Then there are $t_n\geq 0$ and $x_n^*\in\hat\partial\phi(\bar x)$ such that $t_nx^*_n\rightarrow x^*$. We claim that $\{t_n\}$ is bounded. (Otherwise, we can assume that $t_n\rightarrow+\infty$ and thus $x^*_n\rightarrow 0$. This implies that $0\in\hat\partial\phi(\bar x)$, which contradicts $0\not\in\hat\partial\phi(\bar x)$). Since $X$ is finite-dimensional and $\hat\partial\phi(\bar x)$ is bounded, without loss of generalization, we can assume that $t_n\rightarrow t_0\geq 0$ and $x_n^*\rightarrow x_0^*\in\hat\partial\phi(\bar x)$ by the norm-closed property of $\hat\partial\phi(\bar x)$. Then $x^*=t_0x_0^*\in [0,+\infty)\hat\partial\phi(\bar x)$ and thus $[0,+\infty)\hat \partial\phi(\bar x)$ is closed. The proof is complete.\hfill$\Box$

\begin{pro}
Let $\bar x\in {\rm bd}(S)$ and $\tau>0$. Suppose that $X$ is finite dimensional and $\hat\partial\phi(\bar x)$ is bounded. Then inequality \eqref{4-1} satisfies the Fr\'echet $\tau$-strong BCQ at $\bar x$ if and only if
\begin{equation}\label{4-9}
  d\big(h,\overline{\rm co}(T(S,\bar x))\big)\leq\tau\max\left\{0,\sup_{x^*\in\hat\partial\phi(\bar x)}\langle x^*,h\rangle\right\}
\end{equation}
holds for all $h\in X$.

\end{pro}

{\bf Proof.} Noting that $\hat\partial\phi(\bar x)$ is bounded, it follows from Proposition 4.1 that $[0,\tau]\hat \partial\phi(\bar x)$ is closed and convex. This implies that \eqref{4-1} satisfies the Fr\'echet $\tau$-strong BCQ at $\bar x$ if and only if
\begin{equation}\label{4.15}
  \sigma_{\hat N(S,\bar x)\cap B_{X^*}}\leq\sigma_{[0,\tau]\hat \partial\phi(\bar x)}.
\end{equation}
By virtue of \eqref{2.1} again, one can verify that
$$
\sigma_{\hat N(S,\bar x)\cap B_{X^*}}(h)=d\big(h,\overline{\rm co}(T(S,\bar x))\big).
$$
This means that \eqref{4.15} is equivalent to \eqref{4-9}. The proof is complete.\hfill$\Box$

\section{Conclusion}
This paper is devoted to constraint qualifications of the nonconvex inequality defined by a proper lower semicontinuous function. Several types of constraint qualifications involving Clarke/Fr\'echet normal cones and subdifferentials for the nonconvex inequality are investigated. Some necessary and/or sufficient conditions for these constraint qualifications are also provided. When restricted to the convex inequality, these constraint qualifications reduce to BCQ and strong BCQ studied in \cite{H,Z1}. The work in this paper generalizes and extends the study on constraint qualifications from the convex inequality to the nonconvex one. \\

\noindent{\bf Acknowledgment.} The authors wish to thank Professor Xi Yin Zheng for helpful discussions on Remark 3.1 and Fr\'echet BCQ as well as strong BCQ. The authors are also grateful to two anonymous referees for their valuable comments and suggestions which helped us to improve the presentation of this paper.

\end{document}